\documentclass{article}
\usepackage{graphicx}
\usepackage{amssymb,amsmath}
\usepackage{multirow}
\usepackage{bm}
\date{}

\makeatletter \@addtoreset{equation}{section}
\renewcommand{\theequation}{\thesection.\@arabic\c@equation}
\newcommand{\affiliation}[1]{\let\thefootnote\relax\footnote{\mbox{}\\ \noindent {#1}}}

\makeatother
\makeatletter
\def\Ddots{\mathinner{\mkern1mu\raise\p@
 \vbox{\kern7\p@\hbox{.}}\mkern2mu
 \raise4\p@\hbox{.}\mkern2mu\raise7\p@\hbox{.}\mkern1mu}}
\makeatother

\begin{document}
\title{\textbf{On Solving Pentadiagonal Linear Systems via Transformations}}

\author{ A. A. KARAWIA\footnote{ Home Address: Mathematics Department, Faculty of Science, Mansoura University,
Mansoura 35516, Egypt. E-mail:abibka@mans.edu.eg}\\
Computer Science Unit, Deanship of Educational Services, Qassim University,\\
 P.O.Box 6595, Buraidah 51452, Saudi Arabia. \\
              E-mail: kraoieh@qu.edu.sa
}

\maketitle
\begin{abstract}
Many authors studied numerical algorithms for solving the linear systems of the pentadiagonal type. The well-known Fast Pentadiagonal System Solver algorithm is an example of such algorithms. The current article is describes new numerical and symbolic algorithms for solving pentadiagonal linear systems via transformations. New algorithms are natural generalization of the work presented in [Moawwad El-Mikkawy and Faiz Atlan, Algorithms for Solving Linear Systems of Equations of Tridiagonal Type via Transformations, Applied Mathematics, 2014, 5, 413-422]. The symbolic algorithms remove the cases where the numerical algorithms fail. The computational cost of our algorithms is given. Some examples are given in order to illustrate the effectiveness of the proposed algorithms.
All of the experiments are performed on a computer with the aid of programs written in MATLAB.\bigskip
\end{abstract}

\begin{flushleft}\footnotesize
\hspace*{0.9cm}{\textbf{Keywords}:Pentadiagonal matrix; Backward pentadiagonal; Permutation matrix; Linear systems; Algorithm; \hspace*{2.5cm}MATLAB. \\
\hspace*{0.9cm}{\textbf{AMS Subject Classification}:15A15; 15A23; 68W30; 11Y05; 33F10; F.2.1; G.1.0.\\}
}

\end{flushleft}
\newtheorem{alg}{Algorithm}[section]
\newtheorem{cor}{Corollary}[section]

\section{Introduction}

\hspace*{0.5cm}The pentadiagonal linear systems , denoted, by (\textbf{PLS}) take the forms:\\
\begin{equation}
 PX = Y,
 \end{equation}
   where $P$ is $n-by-n$ pentadiagonal matrix given by
\begin{equation}
P=\left(
      \begin{array}{cccccccccc}
        d_1 & a_1 & b_1 & 0 & \ldots & \ldots & \ldots & \ldots & \ldots & 0 \\
        c_2 & d_2 & a_2 & b_2 & 0 & \ldots & \ldots & \ldots & \ldots & 0 \\
        e_3 & c_3 & d_3 & a_3 & b_3 & 0 & \ldots & \ldots & \ldots & 0 \\
        0 & e_4 & c_4 & d_4 & a_4 & b_4 & 0 & \ldots & \ldots & 0 \\
        \vdots & \ddots & \ddots & \ddots & \ddots & \ddots & \ddots & \ddots & \ldots & \vdots \\
        \vdots & \ddots & \ddots & \ddots & \ddots & \ddots & \ddots & \ddots & \ddots & \vdots \\
        \vdots & \ddots & \ddots & \ddots & \ddots & \ddots & \ddots & \ddots & \ddots & \vdots \\
        0 & \ldots & \ldots & \ldots & 0 & e_{n-2} & c_{n-2} & d_{n-2} & a_{n-2}&b_{n-2} \\
        0 & \ldots & \ldots & \ldots & \ldots & 0 & e_{n-1} & c_{n-1} & d_{n-1} & a_{n-1} \\
        0 & \ldots & \ldots & \ldots & \ldots &\ldots & 0 & e_{n} & c_{n} & d_{n}
      \end{array}
    \right),\quad n\ge4.
\end{equation}

and $X =(x_1, x_2,. . ., x_n)^t$, $Y = (y_1, y_2,. . ., y_n)^t$ are vectors of length $n$.\\
\\
This kind of linear systems is well known in the literature [1-11] and often these types of linear systems are widely used in areas of science and engineering, for example in numerical solution of ordinary and partial differential equations (ODE and PDE), interpolation problems, boundary value problems (BVP), parallel computing, Physics, matrix algebra[4-13]. In [7] the author presented an efficient computational algorithm for solving periodic pentadiagonal linear systems. The algorithm is based on the LU factorization of the periodic pentadiagonal matrix. New algorithms are described for solving periodic pentadiagonal linear systems based on the use of any pentadiagonal linear solver and the author described a symbolic algorithm for solving pentadiagonal linear systems[8]. In [9] the authors discussed the general nonsymmetric problem and proposed an algorithm for solving nonsymmetric penta-diagonal Toeplitz linear systems.  A fast algorithm for solving a large system with a symmetric Toeplitz penta-diagonal coefficient matrix is presented[10]. This efficient method is based on the idea of a system perturbation followed by corrections and is competitive with standard methods. In [11] the authors described an efficient computational algorithm and symbolic algorithm for solving nearly pentadiagonal linear systems based on the LU factorization of the nearly pentadiagonal matrix.

In this paper, we show that more efficient algorithms are derived via transformations that can be regarded as a natural generalization of the efficient algorithms in [14].

The current paper is organized as follows: In Section 2, new numerical algorithms for solving a pentadiagonal linear system are presented. New symbolic algorithms for solving a pentadiagonal linear system are constructed in Section 3. In Section 4, three illustrative examples are presented. Conclusions of the work are given in Section 5.\\

\section{Numerical Algorithms for Solving PLS}

\hspace*{0.5cm}In this section we shall focus on the construction of new numerical algorithms for computing the solution of pentadiagonal linear system. For this purpose it is convenient to give five vectors $\boldsymbol{\alpha}=(\alpha_1, \alpha_2, . . ., \alpha_{n-1})$, $\boldsymbol{\beta}=(\beta_1, \beta_2, . . ., \beta_{n-2})$, $\boldsymbol{Z}=(z_1, z_2, . . ., z_n)$, $\boldsymbol{\gamma}=(\gamma_2, \gamma_3, . . ., \gamma_n)$, and $\boldsymbol{\mu}=(\mu_1, \mu_2, . . ., \mu_n)$, where

\begin{equation}
\alpha_i=\left\{\begin{matrix}
\frac{a_1}{\mu_1} & i=1\\
\frac{a_i-\beta_{i-1}\gamma_i}{\mu_i}& \hspace*{1.9cm}i=2,3, ..., n-1,
\end{matrix}\right.
\end{equation}

\begin{equation}
\beta_i=\frac{b_i}{\mu_i}, \hspace*{3.3cm}i=2,3, ..., n-1,
\end{equation}

\begin{equation}
z_i=\left\{\begin{matrix}
\frac{y_1}{\mu_1} & i=1\\
\frac{y_2-z_1\gamma_2}{\mu_2} & i=2\\
\frac{y_i-z_{i-2}e_i-z_{i-1}\gamma_i}{\mu_i}& \hspace*{1.3cm}i=3,4, ..., n,
\end{matrix}\right.
\end{equation}

\begin{equation}
\gamma_i=\left\{\begin{matrix}
c_2 & \hspace*{0.8cm}i=2\\
c_i-\alpha_{i-2}e_i& \hspace*{2cm}i=3,4, ..., n,
\end{matrix}\right.
\end{equation}
and
\begin{equation}
\mu_i=\left\{\begin{matrix}
d_1 & i=1\\
d_2-\alpha_1\gamma_2 & i=2\\
d_i-\beta_{i-2}e_i-\alpha_{i-1}\gamma_i& \hspace*{1.2cm}i=3,4, ..., n.
\end{matrix}\right.
\end{equation}
By using the vectors $\boldsymbol{\alpha}$, $\boldsymbol{\beta}$, $\boldsymbol{Z}$, $\boldsymbol{\gamma}$, and $\boldsymbol{\mu}$, together with the suitable elementary row operations, we see that the system (1.1) may be transformed to the equivalent linear system:\\
\begin{equation}
\left(
      \begin{array}{cccccccccc}
        1 & \alpha_1 & \beta_1 & 0 & \ldots & \ldots & \ldots & \ldots & \ldots & 0 \\
        0 & 1 & \alpha_2 & \beta_2 & 0 & \ldots & \ldots & \ldots & \ldots & 0 \\
        0 & 0 & 1 & \alpha_3 & \beta_3 & 0 & \ldots & \ldots & \ldots & 0 \\
        0 & 0 & 0 & 1 &\alpha_4 & \beta_4 & 0 & \ldots & \ldots & 0 \\
        \vdots & \ddots & \ddots & \ddots & \ddots & \ddots & \ddots & \ddots & \ldots & \vdots \\
        \vdots & \ddots & \ddots & \ddots & \ddots & \ddots & \ddots & \ddots & \ddots & \vdots \\
        \vdots & \ddots & \ddots & \ddots & \ddots & \ddots & \ddots & \ddots & \ddots & \vdots \\
        0 & \ldots & \ldots & \ldots & 0 & 0 & 0 & 1 & \alpha_{n-2}&\beta_{n-2} \\
        0 & \ldots & \ldots & \ldots & \ldots & 0 & 0 & 0 & 1 & \alpha_{n-1} \\
        0 & \ldots & \ldots & \ldots & \ldots &\ldots & 0 & 0 & 0 & 1
      \end{array}
    \right)\left(\begin{array}{c}
             x_1 \\
             x_2 \\
             x_3 \\
             x_4 \\
             \vdots \\
             \vdots \\
             \vdots \\
             x_{n-2} \\
             x_{n-1} \\
             x_n
           \end{array}\right)=\left(\begin{array}{c}
             z_1 \\
             z_2 \\
             z_3 \\
             z_4 \\
             \vdots \\
             \vdots \\
             \vdots \\
             z_{n-2} \\
             z_{n-1} \\
             z_n
           \end{array}\right)
\end{equation}
The transformed system (2.6) is easy to solve by a backward substitution. Consequently, the PLS (1.1) can be solved using the following algorithm:\\

\begin{alg} \textbf{First numerical algorithm for solving pentadiagonal linear system.}\\
\\
\hspace*{1.7cm}To find the solution of PLS (1.1) using the transformed system (2.6), we may proceed as follows:\\
\textbf{INPUT} order of the matrix $n$ and the components $d_i, a_i, b_i, c_i, e_i, f_i,\quad i = 1, 2, . . . , n,(a_n=b_n=b_{n-1}=\hspace*{1.7cm}c_1=e_1=e_2=0)$.\\
\textbf{OUTPUT} The solution vector $x =(x_1, x_2, ..., x_n)^t$.\\
\textbf{Step 1:} Use DETGPENTA algorithm [13] to check the non-singularity of the coefficient matrix of the \hspace*{1.7cm}system (1.3).\\
\textbf{Step 2:} If $det(P)=0$, then Exit and Print Message (''No solutions'') end if.\\
\textbf{Step 3:} Set $\mu_1=d_1$, $\alpha_1=\frac{a_1}{\mu_1}$, $\beta_1=\frac{b_1}{\mu_1}$, and $z_1=\frac{y_1}{\mu_1}$.\\
\textbf{Step 4:} Set $\gamma_2=c_2$, $\mu_2=d_2-\alpha_1\gamma_2$, $\alpha_2=\frac{a_2-\beta_1\gamma_2}{\mu_2}$, $\beta_2=\frac{b_2}{\mu_2}$, and $z_2=\frac{y_2-z_1\gamma_2}{\mu_2}$.\\
\textbf{Step 5:} For i=3,4,...,n-2 do\\
\hspace*{1.7cm}Compute and simplify:\\
    \hspace*{1.7cm}$\gamma_i=c_i-\alpha_{i-2}e_i$,\\
    \hspace*{1.7cm}$\mu_i=d_i-\beta_{i-2}e_i-\alpha_{i-1}\gamma_i$,\\
    \hspace*{1.7cm}$\alpha_i=\frac{a_i-\beta_{i-1}\gamma_i}{\mu_i}$,\\
    \hspace*{1.7cm}$\beta_i=\frac{b_i}{\mu_i}$,\\
    \hspace*{1.7cm}$z_i=\frac{y_i-z_{i-2}e_i-z_{i-1}\gamma_i}{\mu_i}$,\\
    \hspace*{1.4cm}End do.\\
    \hspace*{1.4cm}$\gamma_{n-1}=c_{n-1}-\alpha_{n-3}e_{n-1}$,\\
    \hspace*{1.4cm}$\mu_{n-1}=d_{n-1}-\beta_{n-3}e_{n-1}-\alpha_{n-2}\gamma_{n-1}$,\\
    \hspace*{1.4cm}$\alpha_{n-1}=\frac{a_{n-1}-\beta_{n-2}\gamma_{n-1}}{\mu_{n-1}}$,\\
    \hspace*{1.4cm}$\gamma_n=c_n-\alpha_{n-2}e_n$,\\
    \hspace*{1.4cm}$\mu_n=d_n-\beta_{n-2}e_n-\alpha_{n-1}\gamma_n$,\\
    \hspace*{1.4cm}$z_{n-1}=\frac{y_{n-1}-z_{n-2}e_{n-1}-z_{n-2}\gamma_{n-1}}{\mu_{n-1}}$,\\
    \hspace*{1.4cm}$z_n=\frac{y_n-z_{n-1}e_n-z_{n-1}\gamma_n}{\mu_n}$,\\
\textbf{Step 6:} Compute the solution vector $X=(x_1, x_2, ..., x_n)^t$ using \\
    \hspace*{1.4cm}$x_n=z_n,$ $x_{n-1}=z_{n-1}-\alpha_{n-1}x_n$.\\
\hspace*{1.4cm}For i=n-2, n-3, ...,1 do\\
\hspace*{1.7cm}Compute and simplify:\\
    \hspace*{1.7cm}$x_i=z_i-\alpha_ix_{i+1}-\beta_ix_{i+2}$\\
\hspace*{1.4cm}End do.
\end{alg}
The numerical Algorithm 2.1 will be referred to as \textbf{PTRANS-I} algorithm. The computational cost of \textbf{PTRANS-I} algorithm is $19n-29$ operations. The conditions $\mu_i\ne 0, i=1, 2, ..., n,$ are sufficient for its validity.

In a similar manner, we may consider five vectors $\boldsymbol{\sigma}=(\sigma_2, \sigma_3, . . ., \sigma_n)$, $\boldsymbol{\phi}=(\phi_3, \phi_4, . . ., \phi_n)$, $\boldsymbol{W}=(w_1, w_2, . . ., w_n)$, $\boldsymbol{\rho}=(\rho_1, \rho_2, . . ., \rho_{n-1})$, and $\boldsymbol{\psi}=(\psi_1, \psi_2, . . ., \psi_n)$, where

\begin{equation}
\sigma_i=\left\{\begin{matrix}
\frac{c_n}{\psi_n} & i=n\\
\frac{c_i-\phi_{i+1}\rho_i}{\psi_i}& \hspace*{2.4cm}i= n-1, n-2, ...,2,
\end{matrix}\right.
\end{equation}

\begin{equation}
\phi_i=\frac{e_i}{\psi_i}, \hspace*{3.3cm}i=n,n-1, ..., 3,
\end{equation}

\begin{equation}
w_i=\left\{\begin{matrix}
\frac{y_n}{\psi_n} & i=n\\
\frac{y_{n-1}-w_n\rho_{n-1}}{\psi_{n-1}} & \hspace*{0.6cm}i=n-1\\
\frac{y_i-w_{i+2}b_i-w_{i+1}\rho_i}{\psi_i}& \hspace*{2.4cm}i=n-2, n-3, ...,1,
\end{matrix}\right.
\end{equation}

\begin{equation}
\rho_i=\left\{\begin{matrix}
a_{n-1} & \hspace*{0.8cm}\hspace*{0.4cm}i=n-1\\
a_i-\sigma_{i+2}b_i& \hspace*{2cm}\hspace*{1cm}i=n-2, n-3, ...,1,
\end{matrix}\right.
\end{equation}
and
\begin{equation}
\psi_i=\left\{\begin{matrix}
d_n & i=n\\
d_{n-1}-\sigma_n\rho_{n-1} & \hspace*{0.6cm}i=n-1\\
d_i-\phi_{i+2}b_i-\sigma_{i+1}\rho_i& \hspace*{2.4cm}i=n-2, n-3, ..., 1.
\end{matrix}\right.
\end{equation}
Now we will present another algorithm for solving PLS. As in \textbf{PTRANS-I} algorithm, by using the vectors $\boldsymbol{\sigma}$, $\boldsymbol{\phi}$, $\boldsymbol{W}$, $\boldsymbol{\rho}$, and $\boldsymbol{\psi}$, together with the suitable elementary row operations, we see that the system (1.1) may be transformed to the equivalent linear system:\\

\begin{equation}
\left(
      \begin{array}{cccccccccc}
        1 & 0 & 0 & 0 & \ldots & \ldots & \ldots & \ldots & \ldots & 0 \\
        \sigma_2 & 1 & 0 & 0 & 0 & \ldots & \ldots & \ldots & \ldots & 0 \\
        \phi_3 & \sigma_3 & 1 & 0 & 0 & 0 & \ldots & \ldots & \ldots & 0 \\
        0 & \phi_4 & \sigma_4 & 1 &0 & 0& 0 & \ldots & \ldots & 0 \\
        \vdots & \ddots & \ddots & \ddots & \ddots & \ddots & \ddots & \ddots & \ldots & \vdots \\
        \vdots & \ddots & \ddots & \ddots & \ddots & \ddots & \ddots & \ddots & \ddots & \vdots \\
        \vdots & \ddots & \ddots & \ddots & \ddots & \ddots & \ddots & \ddots & \ddots & \vdots \\
        0 & \ldots & \ldots & \ldots & 0 & \phi_{n-2} & \sigma_{n-2} & 1 & 0&0 \\
        0 & \ldots & \ldots & \ldots & \ldots & 0 & \phi_{n-1} & \sigma_{n-1} & 1 & 0 \\
        0 & \ldots & \ldots & \ldots & \ldots &\ldots & 0 & \phi_n & \sigma_n & 1
      \end{array}
    \right)\left(\begin{array}{c}
             x_1 \\
             x_2 \\
             x_3 \\
             x_4 \\
             \vdots \\
             \vdots \\
             \vdots \\
             x_{n-2} \\
             x_{n-1} \\
             x_n
           \end{array}\right)=\left(\begin{array}{c}
             w_1 \\
             w_2 \\
             w_3 \\
             w_4 \\
             \vdots \\
             \vdots \\
             \vdots \\
             w_{n-2} \\
             w_{n-1} \\
             w_n
           \end{array}\right)
\end{equation}
The transformed system (2.12) is easy to solve by a forward substitution. Consequently, the PLS (1.1) can be solved using the following algorithm:\\

\begin{alg} \textbf{Second numerical algorithm for solving pentadiagonal linear system.}\\
\\
\hspace*{1.5cm}To find the solution of PLS (1.1) using the transformed system (2.12), we may proceed as follows:\\
\textbf{INPUT} order of the matrix $n$ and the components $d_i, a_i, b_i, c_i, e_i, f_i,\quad i = 1, 2, . . . , n,(a_n=b_n=b_{n-1}=\hspace*{1.7cm}c_1=e_1=e_2=0)$.\\
\textbf{OUTPUT} The solution vector $x =(x_1, x_2, ..., x_n)^t$.\\
\textbf{Step 1:} Use DETGPENTA algorithm [13] to check the non-singularity of the coefficient matrix of the \hspace*{1.7cm}system (1.3).\\
\textbf{Step 2:} If $det(P)=0$, then Exit and Print Message (''No solutions'') end if.\\
\textbf{Step 3:} Set $\psi_n=d_n$, $\sigma_n=\frac{c_n}{\psi_n}$, $\phi_n=\frac{e_n}{\psi_n}$, and $w_n=\frac{y_n}{\psi_n}$.\\
\textbf{Step 4:} Set $\rho_{n-1}=a_{n-1}$, $\psi_{n-1}=d_{n-1}-\sigma_n\rho_{n-1}$, $\sigma_{n-1}=\frac{c_{n-1}-\phi_n\rho_{n-1}}{\psi_{n-1}}$, $\phi_{n-1}=\frac{e_{n-1}}{\psi_{n-1}}$, and $w_{n-1}=\hspace*{1.5cm}\frac{y_{n-1}-w_n\rho_{n-1}}{\psi_{n-1}}$.\\
\textbf{Step 5:} For i=n-2, n-3, ...,3 do\\
\hspace*{1.7cm}Compute and simplify:\\
    \hspace*{1.7cm}$\rho_i=a_i-\sigma_{i+2}b_i$,\\
    \hspace*{1.7cm}$\psi_i=d_i-\phi_{i+2}b_i-\sigma_{i+1}\rho_i$,\\
    \hspace*{1.7cm}$\sigma_i=\frac{c_i-\phi_{i+1}\rho_i}{\psi_i}$,\\
    \hspace*{1.7cm}$\phi_i=\frac{e_i}{\psi_i}$,\\
    \hspace*{1.7cm}$w_i=\frac{y_i-w_{i+2}b_i-w_{i+1}\rho_i}{\psi_i}$,\\
    \hspace*{1.4cm}End do.\\
    \hspace*{1.4cm}$\rho_2=a_2-\sigma_4b_2$,\\
    \hspace*{1.4cm}$\psi_2=d_2-\phi_4b_2-\sigma_3\rho_2$,\\
    \hspace*{1.4cm}$\sigma_2=\frac{c_2-\phi_4\rho_2}{\psi_2}$,\\
    \hspace*{1.4cm}$\rho_1=a_1-\sigma_3b_1$,\\
    \hspace*{1.4cm}$\psi_1=d_1-\phi_3b_1-\sigma_2\rho_1$,\\
    \hspace*{1.4cm}$w_2=\frac{y_2-w_4b2-w_3\rho_2}{\psi_2}$,\\
    \hspace*{1.4cm}$w_1=\frac{y_1-w_3b_1-w_2\rho_1}{\psi_1}$,\\
\textbf{Step 6:} Compute the solution vector $X=(x_1, x_2, ..., x_n)^t$ using \\
    \hspace*{1.4cm}$x_=w_1,$ $x_2=w_2-\sigma_2x_1$.\\
\hspace*{1.4cm}For i=3, 4, ...,n do\\
\hspace*{1.7cm}Compute and simplify:\\
    \hspace*{1.7cm}$x_i=w_i-\sigma_ix_{i-1}-\phi_ix_{i-2}$\\
\hspace*{1.4cm}End do.
\end{alg}
The numerical Algorithm 2.2 will be referred to as \textbf{PTRANS-II} algorithm. The computational cost of \textbf{PTRANS-II} algorithm is $19n-29$ operations. Also, the conditions $\psi_i\ne 0, i=1, 2, ..., n,$ are sufficient for its validity.

If $\mu_i= 0$ or $\psi_i= 0$ for any $i\in\{1, 2, ..., n\}$ then \textbf{PTRANS-I} and \textbf{PTRANS-II} algorithm fail to solve pentadiagonal linear systems respectively. So, in the next section, we developed two symbolic algorithms in order to remove the cases where the numerical algorithms fail. The parameter $''p''$ in the following symbolic algorithms is just a symbolic name. It is a dummy argument and its actual value is zero.

\section{Symbolic Algorithms for Solving PLS}

\hspace*{0.5cm}In this section we shall focus on the construction of new symbolic algorithms for computing the solution of pentadiagonal linear systems. The following algorithm is a symbolic version of \textbf{PTRANS-I} algorithm:

\begin{alg} \textbf{First symbolic algorithm for solving pentadiagonal linear system.}\\
\\
\hspace*{1.7cm}To find the solution of PLS (1.1) using the transformed system (2.6), we may proceed as follows:\\
\textbf{INPUT} order of the matrix $n$ and the components $d_i, a_i, b_i, c_i, e_i, f_i,\quad i = 1, 2, . . . , n,(a_n=b_n=b_{n-1}=\hspace*{1.7cm}c_1=e_1=e_2=0)$.\\
\textbf{OUTPUT} The solution vector $x =(x_1, x_2, ..., x_n)^t$.\\
\textbf{Step 1:} Use DETGPENTA algorithm [13] to check the non-singularity of the coefficient matrix of the \hspace*{1.7cm}system (1.3).\\
\textbf{Step 2:} If $det(P)=0$, then Exit and Print Message (''No solutions'') end if.\\
\textbf{Step 3:} Set $\mu_1=d_1$. If $\mu_1=0$ then  $\mu_1=p$ end if.\\
\textbf{Step 4:} Set $\alpha_1=\frac{a_1}{\mu_1}$, $\beta_1=\frac{b_1}{\mu_1}$, $z_1=\frac{y_1}{\mu_1}$ and $\gamma_2=c_2$.\\
\textbf{Step 5:} Set $\mu_2=d_2-\alpha_1\gamma_2$. If $\mu_2=0$ then  $\mu_2=p$ end if.\\
\textbf{Step 6:} Set $\alpha_2=\frac{a_2-\beta_1\gamma_2}{\mu_2}$, $\beta_2=\frac{b_2}{\mu_2}$, and $z_2=\frac{y_2-z_1\gamma_2}{\mu_2}$.\\
\textbf{Step 7:} For i=3,4,...,n-2 do\\
\hspace*{1.7cm}Compute and simplify:\\
    \hspace*{1.7cm}$\gamma_i=c_i-\alpha_{i-2}e_i$,\\
    \hspace*{1.7cm}$\mu_i=d_i-\beta_{i-2}e_i-\alpha_{i-1}\gamma_i$,\\
    \hspace*{1.7cm}If $\mu_i=0$ then  $\mu_i=p$ end if.\\
    \hspace*{1.7cm}$\alpha_i=\frac{a_i-\beta_{i-1}\gamma_i}{\mu_i}$,\\
    \hspace*{1.7cm}$\beta_i=\frac{b_i}{\mu_i}$,\\
    \hspace*{1.7cm}$z_i=\frac{y_i-z_{i-2}e_i-z_{i-1}\gamma_i}{\mu_i}$,\\
    \hspace*{1.4cm}End do.\\
    \hspace*{1.4cm}$\gamma_{n-1}=c_{n-1}-\alpha_{n-3}e_{n-1}$,\\
    \hspace*{1.4cm}$\mu_{n-1}=d_{n-1}-\beta_{n-3}e_{n-1}-\alpha_{n-2}\gamma_{n-1}$. If $\mu_{n-1}=0$ then  $\mu_{n-1}=p$ end if.\\
    \hspace*{1.4cm}$\alpha_{n-1}=\frac{a_{n-1}-\beta_{n-2}\gamma_{n-1}}{\mu_{n-1}}$,\\
    \hspace*{1.4cm}$\gamma_n=c_n-\alpha_{n-2}e_n$,\\
    \hspace*{1.4cm}$\mu_n=d_n-\beta_{n-2}e_n-\alpha_{n-1}\gamma_n$. If $\mu_n=0$ then  $\mu_n=p$ end if. \\
    \hspace*{1.4cm}$z_{n-1}=\frac{y_{n-1}-z_{n-2}e_{n-1}-z_{n-2}\gamma_{n-1}}{\mu_{n-1}}$,\\
    \hspace*{1.4cm}$z_n=\frac{y_n-z_{n-1}e_n-z_{n-1}\gamma_n}{\mu_n}$,\\
\textbf{Step 8:} Compute the solution vector $X=(x_1, x_2, ..., x_n)^t$ using \\
    \hspace*{1.4cm}$x_n=z_n,$ $x_{n-1}=z_{n-1}-\alpha_{n-1}x_n$.\\
\hspace*{1.4cm}For i=n-2, n-3, ...,1 do\\
\hspace*{1.7cm}Compute and simplify:\\
    \hspace*{1.7cm}$x_i=z_i-\alpha_ix_{i+1}-\beta_ix_{i+2}$\\
\hspace*{1.4cm}End do.\\
\textbf{Step 9:} Substitute $p=0$ in all expressions of the solution vector $x_i, i=1, 2, ...,n$.
\end{alg}
The symbolic Algorithm 3.1 will be referred to as \textbf{SPTRANS-I} algorithm.\\
Now we are going to give the symbolic version of \textbf{PTRANS-II} algorithm:

\begin{alg} \textbf{Second symbolic algorithm for solving pentadiagonal linear system.}\\
\\
\hspace*{1.5cm}To find the solution of PLS (1.1) using the transformed system (2.12), we may proceed as follows:\\
\textbf{INPUT} order of the matrix $n$ and the components $d_i, a_i, b_i, c_i, e_i, f_i,\quad i = 1, 2, . . . , n,(a_n=b_n=b_{n-1}=\hspace*{1.7cm}c_1=e_1=e_2=0)$.\\
\textbf{OUTPUT} The solution vector $x =(x_1, x_2, ..., x_n)^t$.\\
\textbf{Step 1:} Use DETGPENTA algorithm [13] to check the non-singularity of the coefficient matrix of the \hspace*{1.7cm}system (1.3).\\
\textbf{Step 2:} If $det(P)=0$, then Exit and Print Message (''No solutions'') end if.\\
\textbf{Step 3:} Set $\psi_n=d_n$. If $\psi_n=0$ then  $\psi_n=p$ end if.\\
\textbf{Step 4:} $\sigma_n=\frac{c_n}{\psi_n}$, $\phi_n=\frac{e_n}{\psi_n}$, $w_n=\frac{y_n}{\psi_n}$ and $\rho_{n-1}=a_{n-1}$.\\
\textbf{Step 5:} Set $\psi_{n-1}=d_{n-1}-\sigma_n\rho_{n-1}$.If $\psi_{n-1}=0$ then  $\psi_{n-1}=p$ end if.\\
\textbf{Step 6:} $\sigma_{n-1}=\frac{c_{n-1}-\phi_n\rho_{n-1}}{\psi_{n-1}}$, $\phi_{n-1}=\frac{e_{n-1}}{\psi_{n-1}}$, and $w_{n-1}=\frac{y_{n-1}-w_n\rho_{n-1}}{\psi_{n-1}}$.\\
\textbf{Step 7:} For i=n-2, n-3, ...,3 do\\
\hspace*{1.7cm}Compute and simplify:\\
    \hspace*{1.7cm}$\rho_i=a_i-\sigma_{i+2}b_i$,\\
    \hspace*{1.7cm}$\psi_i=d_i-\phi_{i+2}b_i-\sigma_{i+1}\rho_i$,\\
    \hspace*{1.7cm}If $\psi_i=0$ then  $\psi_i=p$ end if.\\
    \hspace*{1.7cm}$\sigma_i=\frac{c_i-\phi_{i+1}\rho_i}{\psi_i}$,\\
    \hspace*{1.7cm}$\phi_i=\frac{e_i}{\psi_i}$,\\
    \hspace*{1.7cm}$w_i=\frac{y_i-w_{i+2}b_i-w_{i+1}\rho_i}{\psi_i}$,\\
    \hspace*{1.4cm}End do.\\
    \hspace*{1.4cm}$\rho_2=a_2-\sigma_4b_2$,\\
    \hspace*{1.4cm}$\psi_2=d_2-\phi_4b_2-\sigma_3\rho_2$. If $\psi_2=0$ then  $\psi_2=p$ end if.\\
    \hspace*{1.4cm}$\sigma_2=\frac{c_2-\phi_4\rho_2}{\psi_2}$,\\
    \hspace*{1.4cm}$\rho_1=a_1-\sigma_3b_1$,\\
    \hspace*{1.4cm}$\psi_1=d_1-\phi_3b_1-\sigma_2\rho_1$. If $\psi_1=0$ then  $\psi_i=p$ end if.\\\\
    \hspace*{1.4cm}$w_2=\frac{y_2-w_4b2-w_3\rho_2}{\psi_2}$,\\
    \hspace*{1.4cm}$w_1=\frac{y_1-w_3b_1-w_2\rho_1}{\psi_1}$,\\
\textbf{Step 8:} Compute the solution vector $X=(x_1, x_2, ..., x_n)^t$ using \\
    \hspace*{1.4cm}$x_=w_1,$ $x_2=w_2-\sigma_2x_1$.\\
\hspace*{1.4cm}For i=3, 4, ...,n do\\
\hspace*{1.7cm}Compute and simplify:\\
    \hspace*{1.7cm}$x_i=w_i-\sigma_ix_{i-1}-\phi_ix_{i-2}$\\
\hspace*{1.4cm}End do.\\
\textbf{Step 9:} Substitute $p=0$ in all expressions of the solution vector $x_i, i=1, 2, ...,n$.
\end{alg}
The symbolic Algorithm 3.2 will be referred to as \textbf{SPTRANS-II} algorithm.\\
\begin{cor}(generalization version of Corollary 2.1 in [14]) Let $\hat{P}$ be the backward matrix of the pentadiagonal matrix $P$ in (1.2), and given by:
\begin{equation}
\hat{P}=\left(
      \begin{array}{cccccccccc}
        0 & \ldots & \ldots & \ldots & \ldots & \ldots & 0 & b_1 & a_1 & d_1 \\
        0 & \ldots & \ldots & \ldots & \ldots & 0 & b_2 & a_2 & d_2 & c_2 \\
        0 & \ldots & \ldots & \ldots & 0 & b_3 & a_3 & d_3 & c_3 & e_3 \\
        0 & \ldots & \ldots & 0& b_4 & a_4 & d_4 & c_4 & e_4 & 0 \\
        \vdots & \ldots & \Ddots & \Ddots & \Ddots & \Ddots & \Ddots & \Ddots & \Ddots & \vdots \\
        \vdots & \Ddots & \Ddots & \Ddots & \Ddots & \Ddots & \Ddots & \Ddots & \ldots & \vdots \\
        \vdots & \Ddots & \Ddots & \Ddots & \Ddots & \Ddots & \Ddots & \ldots & \ldots & \vdots \\
        b_{n-2} & a_{n-2} & d_{n-2} & c_{n-2} & e_{n-2} & 0 & \ldots & \ldots & \ldots &0 \\
        a_{n-1} & d_{n-1} & c_{n-1} & e_{n-1} & 0 & \ldots & \ldots & \ldots & \ldots & 0 \\
        d_n & c_n & e_n & 0 & \ldots &\ldots & \ldots& \ldots & \ldots & 0
      \end{array}
    \right),\quad n\ge4.
\end{equation}
Then the backward pentadiagonal linear system
\begin{equation}
\hat{P}V=Y,\quad V=(v_1, v_2, ..., v_n)^t.
\end{equation}
has the solution: $v_i=x_{n-i+1}, i=1,2, ..., \lfloor n\rfloor$, where $\lfloor j\rfloor$ is the floor function of $j$ and $X=(x_1,x_2, ..., x_n)^t$ is the solution vector of the linear system (1.1).\\
\textbf{Proof:} Consider the $n\times n$ permutation matrix $M$ defined by:
\begin{equation}
M=\left(\begin{array}{ccccc}
    0 & \cdots & \cdots & 0 & 1 \\
    \vdots &  &  & 1 & 0 \\
    \vdots &  & \Ddots &  & \vdots \\
    0 & 1 &  &  & \vdots \\
    1 & 0 & \cdots & \cdots & 0
  \end{array}\right)
\end{equation}
For this matrix, we have:
\begin{equation}
M^{-1}=M.
\end{equation}
Since
\begin{equation}
\hat{P}=PM
\end{equation}
Then using (3.4) and (3.5), the result follows.
\end{cor}
\begin{cor}(generalization version of Corollary 2.2 in [14]) The determinants of the coefficient matrices $P$ and $\hat{P}$ in (1.2) and (3.1) are given respectively by:
\begin{equation}
det(P)=\prod_{i=1}^n\mu_i=\prod_{i=1}^n\psi_i
\end{equation}
and
\begin{equation}
det(\hat{P})=(-1)^{\frac{n(n-1)}{2}}\prod_{i=1}^n\mu_i=(-1)^{\frac{n(n-1)}{2}}\prod_{i=1}^n\psi_i
\end{equation}
where $\mu_1, \mu_2, ..., \mu_n$ and $\psi_1, \psi_2, ..., \psi_n$ satisfy (2.5) and (2.11) respectively.\\
\\
\textbf{Proof:} Using (2.6), (2.12) and (3.5), the result follows.
\end{cor}
\section{ILLUSTRATIVE EXAMPLES}
       In this section we are going to give three examples for the sake of illustration. All experiments performed
 in MATLAB R2014a with an Intel(R) Core(TM) i7-4700MQ CPU@2.40GHz 2.40 GHz.\\
\\
\textbf{Example 4.1.(Case I: $\mu_i\ne 0$ and $\psi_i\ne 0$ for all $i$)}\\ Find the solution of the following pentadiagonal linear system of size 10

\begin{equation}
\left(
  \begin{array}{cccccccccc}
    1 & 2 & 1 & 0 & 0 & 0 & 0& 0 & 0 & 0 \\
    3 & 2 & 2 & 5 & 0 & 0 & 0& 0 & 0 & 0 \\
    1 & 2 & 3 & 1 & -2 & 0 & 0& 0 & 0 & 0 \\
    0 & 3 & 1 & -4 & 5 & 1 & 0 & 0 & 0 & 0\\
    0 & 0 & 1 & 2 & 5 & -7 & 5 & 0 & 0 & 0\\
    0 & 0 & 0 & 5 & 1 & 6 & 3 & 2 & 0 & 0\\
    0 & 0 & 0 & 0 & 2 & 2 & 7 & -1 & 4 & 0\\
    0 & 0 & 0 & 0 & 0 & 2 & 1 & -1 & 4 & -3\\
    0 & 0 & 0 & 0 & 0 & 0 & 2 & -2 & 1 & 5\\
    0 & 0 & 0 & 0 & 0 & 0 & 0 & -1 & 4 & 8
  \end{array}
\right)\left(
         \begin{array}{c}
           x_1 \\
           x_2 \\
           x_3 \\
           x_4 \\
           x_5 \\
           x_6 \\
           x_7 \\
           x_8\\
           x_9\\
           x_{10}
         \end{array}
       \right)=\left(
                 \begin{array}{c}
                   8 \\
                   33 \\
                   8 \\
                   24 \\
                   29 \\
                   82 \\
                   71\\
                   17\\
                   57\\
                   108
                 \end{array}
               \right)
\end{equation}

\textbf{Solution:} We have\\
$n=10$, $d=(1, 2, 3, -4, 5, 6, 7, -1, 1, 8)^t$, $a=(2,2,1,5,-7,3,-1,4,5)^t$, $b=(1,5,-2,1,5,2,4,-3)^t$, $c=(0,3,2,1,2,1,2,1,-2,4)^t$,$e=(0,0,1,3,1,5,2,2,2,-1)^t$, and $y=(8,33,8,24,29,98,99,17,57,108)^t$.\\
i)- Applying the \textbf{PTRANS-I} algorithm, it yields
\begin{itemize}
  \item { $\mathbf{\mu}=( 1,-4,2,-\frac{3}{8},27,\frac{245}{9},\frac{3289}{441},-\frac{335}{383},-\frac{2897}{484},\frac{3439}{279} )^t$, $det(P)=\prod_{i=1}^{10}\mu_i={\frac {4989610795975}{4701708}}$.}
  \item {\textbf{PTRANS-I}(n,d,a,b,c,e,y)=$(1,2,3,4,5,6,7,8,9,10)^t$.}
\end{itemize}\bigskip
ii)- Applying the \textbf{PTRANS-II} algorithm, it yields
\begin{itemize}
  \item { $\mathbf{\psi}=(-{\frac {6213}{3613}},{\frac {1603}{1405}},{\frac {1487}{433}},-{
\frac {5173}{239}},{\frac {383}{156}},{\frac {988}{161}},{\frac {69}{
11}},-{\frac {77}{12}},-\frac{3}{2},8)^t$, $det(P)=\prod_{i=1}^{10}\psi_i={\frac {557494642026514353}{525327436055}}$.}
  \item {\textbf{PTRANS-II}(n,d,a,b,c,e,y)=$(1,2,3,4,5,6,7,8,9,10)^t$.}
\end{itemize}\bigskip
\textbf{Example 4.2.}(Case II: $\mu_i= 0$ and $\psi_i= 0$ for some $i$)\\ Find the solution of the following pentadiagonal linear system of size 4

\begin{equation}
\left(\begin{array}{cccc} 3 & 2 & 1 & 0\\ -3 & -2 & 7 & 1\\ 3 & 2 & -1 & 5\\ 0 & 1 & 2 & 3 \end{array}\right)\left(
         \begin{array}{c}
           x_1 \\
           x_2 \\
           x_3 \\
           x_4
         \end{array}
       \right)=\left(
                 \begin{array}{c}
                   6 \\
                   3 \\
                   9 \\
                   6
                 \end{array}
               \right)
\end{equation}
\textbf{Solution:} We have\\
$n=4$, $d=(3,-2,-1,3)^t$, $a=(2,7,5)^t$, $b=(1,1)^t$, $c=(0,-3,2,2)^t$,\\$e=(0,0,3,1)^t$, and $y=(6,3,9,6)^t$.\\
\\
The numerical algorithms \textbf{PTRANS-I} and \textbf{PTRANS-II} fail to solve the pentadiagonal linear system (4.2) since $\mu_2=0$.\\
\\
i)- Applying the \textbf{SPTRANS-I} algorithm, it yields
\begin{itemize}
  \item {$\mathbf{\mu}=(3, p, -2,\frac{8p - 21}{p})^t$. $det(P)=(\prod_{i=1}^4\mu_i)_{p=0}=126$. }
  \item {\textbf{SPTRANS-I}(n,d,a,b,c,e,y)=$((\frac{(25p - 42)}{(16p - 42)}, \frac{-21}{(8p - 21)}, \frac{21(p - 2)}{2(8p - 21)},\frac{(9p - 21)}{(8p -21)})^t)_{p=0}=(1, 1, 1, 1)^t$.}
\end{itemize}
ii)- Applying the \textbf{SPTRANS-II} algorithm, it yields
\begin{itemize}
  \item {$\mathbf{\psi}=(\frac{21}{4}, \frac{-24}{13}, \frac{-13}{3}, 3)^t$. $det(P)=\prod_{i=1}^4\psi_i=126$. }
  \item {\textbf{SPTRANS-II}(n,d,a,b,c,e,y)=$(1, 1, 1, 1)^t$.}
\end{itemize}
\textbf{Example 4.3.} We consider the following $n\times n$ pentadiagonal linear system in order to demonstrate the efficiency of algorithms 3.1 and 3.2.\\
\\
$
\left(
  \begin{array}{cccccccccc}
    9 & -4 & 1 & 0 & \cdots & \cdots & \cdots & \cdots & \cdots & 0 \\
    -4 & 6 & -4 & 1 & 0& &  &  &  & 0 \\
    1 & -4 & 6 & -4 & 1 & 0 &  &  &  & 0 \\
    0& 1 & -4 & 6 & -4 & 1 & \ddots &  &   & 0 \\
    \vdots & \ddots& \ddots & \ddots & \ddots & \ddots &\ddots  & \ddots &  & \vdots \\
     \vdots&  & \ddots & \ddots & \ddots & \ddots & \ddots & \ddots & \ddots & \vdots \\
    \vdots &  &  & \ddots & \ddots & \ddots & \ddots & \ddots & \ddots & 0 \\
    \vdots &  &  & &0 &1 & -4 & 6 & -4 & 1 \\
    \vdots &  &  &  &  & 0 & 1 & -4 & 5 & -2 \\
    0 & \cdots & \cdots & \cdots &\cdots  &\cdots & 0 & 1 & -2 & 1 \\\\
  \end{array}
\right)\left(
          \begin{array}{c}
            x_1 \\
            x_2 \\
            x_3 \\
            x_4 \\
            \vdots \\
            \vdots \\
            \vdots\\
            x_{n-2} \\
            x_{n-1} \\
            x_n \\
          \end{array}
        \right)=\left(
          \begin{array}{c}
            6 \\
            -1\\
            0 \\
            0 \\
            \vdots \\
            \vdots\\
            \vdots \\
            0 \\
            0 \\
            0 \\
          \end{array}
        \right)\\
        $

It can be verified that the exact solution is $\mathbf{x}=(1, 1, \ldots, 1)^t$. In Table 1 we give a comparison of the absolute error and running time between \textbf{PTRANS-I}, \textbf{PTRANS-II} algorithms, \textbf{Algorithm 3}[8], and "A$\backslash$b" function in Matlab to compute $\bar{x}$. for different sizes.\\
\newpage
$$Table1.$$
$$
\begin{tabular}{|c|c|c|c|c|}
  \hline

  \multirow{2}{*}{n} & \multicolumn{4}{c|}{$||\mathbf{x}-\bar{x}||_\infty$ and CPU time(S)} \\ \cline{2-5}
   & \textbf{PTRANS-I} & \textbf{PTRANS-II} & \textbf{Algorithm 3}[8] & A$\backslash$b(MATLAB) \\ \hline
  500 & $1.5856\times 10^{-7}\quad  0.0069$ & $0\hspace*{1.2cm}  0.0086$ & $6.8579\times 10^{-8}\quad  0.0048$ & $9.98\times 10^{-8}\quad 0.0023$\\ \hline
  5000 & $8.3674\times 10^{-4}\quad  0.0062$ & $0\hspace*{1.2cm}  0.0391$ & $3.0253\times 10^{-4}\quad  0.0057$ &$2.50\times 10^{-4}\quad 0.7548$\\ \hline
  10000 & $0.0058\hspace*{1.2cm}  0.0114$ & $0\hspace*{1.2cm} 0.0511$ & $0.0052\hspace*{1.2cm}  0.0101$ & $0.0106\hspace*{1.2cm} 4.5464$\\ \hline
  50000 & $2.1415\hspace*{1.2cm} 0.0308$ & $0\hspace*{1.2cm} 0.1687$ & $7.9056\hspace*{1.2cm}  0.0119$ & $0.0159\hspace*{1.2cm} 655.51$\\ \hline
  \end{tabular}
$$
Table 1 indicates that the value of absolute error for large values of $n$ for $\textbf{PTRANS-II}$ is less than the other algorithms. Also note that from Table 1 it is clear that the value of the runing time for $\textbf{Algorithm 3}[8]$ is less than the other algorithms for large sizes.\\

\textbf{Acknowledgements} The author wishes to thank anonymous referees for useful comments that enhanced the quality of this paper. The author is grateful to Prof. Dr. M.E.A. El-Mikkawy for providing him with the references [14] and [15].

\section{CONCLUSION}

There are many numerical algorithms in current use for solving linear systems of pentadiagonal type. All numerical algorithms including the PTRANS-I and PTRANS-II algorithms of the current paper, fail to solve the pentadiagonal linear system if $\mu_i=0$ and $\psi_i=0$ for any $i\in\{1,2,...,n\}$. The symbolic algorithms SPTRANS-I and SPTRANS-II of the current paper are constructed in order to remove the cases where the numerical algorithms fail. From some numerical examples we have learned that SPTRANS-II algorithm works as well as $\textbf{Algorithm 3}[8]$ and (A$\backslash$y)MATLAB algorithms. Hence, it may become a useful tool for solving linear systems of pentadiagonal type.

\end{document}